\documentclass[a4paper,12pt]{amsart}

\usepackage{amsthm}

\usepackage{amssymb} % symboler fra AMS
\usepackage{latexsym} % LaTeX symboler
\usepackage{amsfonts} % fonte fra AMS
\usepackage{amsmath} % matematik fonte fra AMS
\usepackage{eucal} % 'krlle'-fonte (skriv \mathcal)
\usepackage{bm} % fed skrift i matematik
\usepackage{bbm} % skrift med dobbelt streg til f.eks. talomr�er
\usepackage{graphicx} % mulighed for at inkludere grafik
\usepackage[english]{varioref} % smarte referencer - se nederst
\usepackage[nice]{nicefrac} % p�ere brker
\usepackage[all]{xy}

\newcommand{\Ad}{\text {\rm Ad}}
\newcommand{\Hom}{\text {\rm Hom}}

\newcommand{\supp}{\text {\rm supp}}

\newcommand{\Cyc}{\text{\rm Cyc}}

\newcommand{\kk}{\mathbf k}

\newcommand{\cb}{\mathcal B}

\newcommand{\co}{\mathcal O}

\newcommand{\NN}{\mathbb N}
\newcommand{\ZZ}{\mathbb Z}

\def\ge{\geqslant}
\def\le{\leqslant}
\def\a{\alpha}

\def\d{\delta}

\def\e{\epsilon}

\def\o{\omega}

\def\l{\lambda}

\def\i{^{-1}}

\theoremstyle{plain}%default
\newtheorem*{thm*}{Theorem} % thm* er nu uden nummer

 \newtheorem*{rmk}{Remark}
 
 \newtheorem*{clm}{Claim}

\newtheorem*{th1}{Theorem 1.3}
\newtheorem*{th2}{Theorem 1.5}
\newtheorem*{th3}{Theorem 1.7}
\newtheorem*{th4}{Proposition 1.9}
\newtheorem*{th5}{Lemma 1.11}

\title{On the affineness of Deligne-Lusztig varieties}
\author{Xuhua He}
\address{Department of Mathematics, Stony Brook University, Stony Brook,
NY 11794, USA}
\email{hugo@math.sunysb.edu}
\thanks{The author is partially supported by NSF grant DMS-0700589}

\subjclass[2000]{20C33, 20F55}

\begin{document}

\begin{abstract}
We prove that the Deligne-Lusztig variety associated to minimal
length elements in any $\d$-conjugacy class of the Weyl group is
affine, which was conjectured by Orlik and Rapoport in \cite{OR}.
\end{abstract}
\maketitle

\subsection*{1.1 Notations} Let $\kk$ be an algebraic closure of the finite prime
field $\mathbb F_p$ and $G$ be a connected reductive algebraic
group over $\kk$ with an endomorphism $F: G \to G$ such that some
power $F^d$ of $F$ is the Frobenius endomorphism relative to a
rational structure over a finite subfield $\kk_0$ of $\kk$. Let $q$
be the positive number with $q^d=|\kk_0|$.

We fix a $F$-stable Borel subgroup $B$ and a $F$-stable maximal
torus $T \subset B$. Let $\Phi$ be the set of roots and $(\a_i)_{i
\in I}$ be the set of simple roots corresponding to $(B, T)$. For $i
\in I$, let $\o^\vee_i$ be the corresponding fundamental coweight.
Let $W=N(T)/T$ be the Weyl group and $(s_i)_{i \in I}$ be the set of
simple reflections. For $w \in W$, let $l(w)$ be the length of $w$.
Since $(B, T)$ is $F$-stable, $F$ induces a bijection on $I$ and an
automorphism on $W$. We denote the induced maps on $W$ and $I$ by
$\d$. Now $\d$ also induces isomorphisms on the set of characters
$X=\Hom(T, G_m)$ and the set of cocharacters $X^\vee=\Hom(G_m, T)$
which we also denote by $\d$. Then it is easy to see that $F^* \mu=q
\cdot \d \i(\mu)$ for $\mu \in X^\vee$.

For $J \subset I$, let $\Phi_J$ be the set of roots generated by $\{\a_j\}_{j \in J}$ and $W_J$ be the subgroup of $W$ generated by $\{s_j\}_{j \in J}$. Let $W^J$ be the set of minimal length coset representatives for $W/W_J$. The unique maximal element in $W$ will be denoted by $w_0$ and the unique maximal element in $W_J$ will be denoted by $w_0^J$.

Let $\a_0=\sum_{i \in I} n_i \a_i$ be the highest root and
$n_0=\sum_{i \in I} n_i$.

\subsection*{1.2} Let $\cb$ be the set of Borel subgroups of $G$. For $w
\in W$, let $O(w)=\{({}^g B, {}^{g \dot w} B); g \in G\}$ be the
$G$-orbit on $\cb \times \cb$ that corresponding to $w$. Set
\[X(w)=\{ B' \in \cb; (B', F(B')) \in O(w)\}.\] This is the
Deligne-Lusztig variety associated to $w$ (see \cite[1.4]{DL}). It
is known that $X(w)$ is a variety of pure dimension $l(w)$ (see {\it
loc.cit.}) and is quasi-affine (see \cite{Ha}). It is also known
that when $q \ge h$ (where $h$ is the Coxeter number), then $X(w)$
is affine (see \cite[Theorem 9.7]{DL}).

The main result we will prove in this note is the following

\begin{th1} Let $w \in W$ be a minimal length element in the
$\d$-conjugacy class $\{x w \d(x) \i; x \in W\}$. Then $X(w)$ is
affine.
\end{th1}

\begin{rmk}
The case where $w$ is a Coxeter element was proved by Lusztig in \cite[Corollary 2.8]{L} in a geometric way and the cases for split classical groups were proved  by Orlik and Rapoport in
\cite[section 5]{OR} by finding a minimal length element in each $\d$-conjugacy that satisfies the criterion \cite[Theorem 9.7]{DL}. Our approach is motivated by the approach of Orlik and Rapoport. However, a main difference is the way of choosing the minimal length elements. We will discuss it in more detail in 1.14.
\end{rmk}

Before discussing the proof of the theorem above, we first recall
some results on the minimal length elements.

\subsection*{1.4} We follow the notations in \cite[section 3.2]{GP2}.

Let $w, w' \in W$ and $j \in I$, we write $w \xrightarrow{s_j}_{\d} w'$
if $w'=s_j w \d(s_j)$ and $l(w') \le l(w)$. If $w=w_0, w_1, \cdots,
w_n=w'$ is a sequence of elements in $W$ such that for all $k$, we
have $w_{k-1} \xrightarrow{s_j}_{\d} w_k$ for some $j \in I$, then we
write $w \rightarrow_{\d} w'$.

We call $w, w' \in W$ {\it elementarily strongly $\d$-conjugate} if
$l(w)=l(w')$ and there exists $x \in W$ such that $w'=x w \d(x)
\i$ and $l(x w)=l(x)+l(w)$ or $l(w \d(x) \i)=l(x)+l(w)$. We
call $w, w'$ {\it strongly $\d$-conjugate} if there is a sequence
$w=w_0, w_1, \cdots, w_n=w'$ such that $w_{i-1}$ is elementarily
strongly $\d$-conjugate to $w_i$. We will write $w \sim_{\d} w'$ if $w$
and $w'$ are strongly $\d$-conjugate.

If $w \sim_{\d} w'$ and $w \rightarrow_{\d} w'$, then we say that $w$ and $w'$ are in the same $\d$-cyclic shift class and write $w \approx_{\d} w'$. For $w \in W$, set $$\Cyc_{\d}(w)=\{w' \in W; w \approx_{\d} w'\}.$$

The following result was proved in \cite[Theorem 1.1]{GP1} for the
usual conjugacy classes and in \cite[Theorem 2.6]{GKP} for the
twisted conjugacy classes.

\begin{th2} Let $\co$ be a $\d$-conjugacy class in $W$ and $\co_{\min}$ be the set of minimal length
elements in $\co$. Then

(1) For each $w \in \co$, there exists $w' \in \co_{\min}$ such that
$w \rightarrow_{\d} w'$.

(2) Let $w, w' \in \co_{\min}$, then $w \sim_{\d} w'$.
\end{th2}

\subsection*{1.6} In general, $\co_{\min}$  might be a union of several $\d$-cyclic shift classes. However, for some special $\d$-conjugacy classes, we have a better result. Let us first introduce some notations.

For $w \in W$, set $\supp_{\d}(w)=\cup_{n \ge 0} \d^n \supp(w)$. Then $\supp_{\d}(w)$ is the minimal $\d$-stable subset of $I$ such that $w \in W_{\supp_{\d}(w)}$.

A $\d$-conjugacy class $\co$ of $W$ is called {\it cuspidal} if $\co \cap W_J=\varnothing$ for all proper $\d$-stable subset $J$ of $I$.

The following result was proved in \cite[Theorem 3.2.7]{GP2} for the usual conjugacy classes, in \cite[section 6]{GKP} for twisted conjugacy classes of exceptional groups and in \cite[Theorem 7.5]{He} for twisted conjugacy classes of classical groups.

\begin{th3} Let $\co$ be a $\d$-conjugacy class and $w \in \co_{\min}$. Then 

(1) If $\supp_{\d}(w)=I$, then $\co$ is cuspidal.

(2) If $\co$ is cuspidal, then $\co_{\min}=\Cyc_{\d}(w)$.
\end{th3}

\subsection*{1.8} By \cite[Lemma 2.6]{OR}, if $w \approx_{\d} w'$, then  $X(w)$ and
$X(w')$ are universally homeomorphic. In particular, if $X(w)$ is
affine, then $X(w')$ is also affine for any $w' \approx_{\d} w$. (However, it is unknown if the same result holds when $w' \sim_{\d} w$.)

By \cite[Theorem 9.7]{DL}, to prove our main theorem, it suffices to prove the following result.

\begin{th4} Let $C=\{ \mu \in X^\vee
\otimes \mathbb R; \a_i(\mu)>0 \text{ for } i \in I\}$ be the
fundamental chamber corresponding to $B$. Then for any $\d$-conjugacy 
class $\co$ of $W$ and $w' \in \co_{\min}$, there exists $w \approx_{\d} w'$ and $\mu \in X^\vee \otimes \mathbb R$ such that $\a(\mu)>0$ for $\a>0$ with $w \a<0$ and $F^* \mu-w \cdot \mu \in C$.
\end{th4}

\subsection*{1.10 Reduction to cuspidal classes} 
Assume that the Proposition 1.9 holds in the case where $\co$ is
cuspidal. We will prove now that it holds in general.

Let $\co$ be an $\d$-conjugacy class of $W$ and $w' \in \co_{\min}$. Let $J=\supp_{\d}(w')$ and $\co'$ be the $\d$-conjugacy class of $W_J$ that contains $w'$. Then $w' \in \co'_{\min}$ and $\co'_{\min} \subset \co_{\min}$. By Theorem 1.7 (1), $\co'$ is a cuspidal $\d$-conjugacy class of $W_J$. Notice that if $x \in W_J$ and $\a \in \Phi$ with $\a>0$ and $x \a<0$, then $\a \in \Phi_J$. By our hypothesis, there exist $w \in \co'_{\min}$ and
$\mu=\sum_{i \in J} m_i \o^\vee_i$ for some $m_i \in \mathbb R$ such
that $\a(\mu)>0$ for $\a>0$ with $w \a<0$ and $\a_i \bigl(F^* \mu-w
\cdot \mu \bigr)>0$ for all $i \in J$. Set $\l=\mu+m \sum_{i \notin
J} \o^\vee_i$ for $m \gg 0$. Then $\a(\l)=\a(\mu)$ for $\a>0$ with
$w \a<0$. 

If $i \in J$, then $\a_i \bigl(F^* \l-w \cdot \l
\bigr)=\a_i \bigl(F^* \mu-w \cdot \mu \bigr)>0$. 

If $i \notin J$,
then $w \i (\a_i)=\a_i+\sum_{j \in J} a_j \a_j$ for some $a_j \in \NN$ with $\sum_{j \in J} a_j \le n_0$. Hence $$\a_i \bigl(F^* \l-w \cdot \l \bigr)=q m-w \i (\a_i)(\l)>(q-1)
m-n_0 \max_{i \in J} |m_i|>0.$$

Therefore, to prove Proposition 1.9, it suffices to prove the following statement.

\begin{th5} For any cuspdial $\d$-conjugacy 
class $\co$ of $W$, there exists $w \in \co_{\min}$ and $\mu \in X^\vee \otimes \mathbb R$ such that $\a(\mu)>0$ for $\a>0$ with $w \a<0$ and $F^* \mu-w \cdot \mu \in C$.
\end{th5}

\subsection*{1.12 Reduction to irreducible types} Assume that Lemma 1.11 
holds in the case where $(\Phi, I)$ is irreducible.
We will prove now that it holds in general.

{\it Step 1}. Assume that $\Phi=\Phi_1 \sqcup \Phi_2 \sqcup \cdots
\sqcup \Phi_r$, where each $\Phi_i$ is irreducible and generated by
$I_i=I \cap \Phi_i$ such that $\d(I_i)=I_{i+1}$ for $i<r$ and
$\d(I_r)=I_1$. Then $W=W_1 \times W_2 \times \cdots \times W_r$ and
we may regard $W_1$ as a subgroup of $W$ in the natural way. In this
case, $\co \cap W_1$ is a cuspidal $\d^r$-conjugacy class of $W_1$. Notice that if $x \in W_1$
and $\a \in \Phi$ with $\a>0$ and $x \a<0$, then $\a \in \Phi_1$. By
assumption, there exists $w \in \co_{\min} \cap W_1=(\co \cap
W_1)_{\min}$ and $\mu=\sum_{i \in I_1} m_i \o^\vee_i$ for some $m_i
\in \mathbb R$ such that $\a(\mu)>0$ for $\a>0$ with $w \a<0$ and
$\a_i \bigl((F^*)^r \mu-w \cdot \mu \bigr)>0$ for all $i \in I_1$.
Therefore, there exists $\e_i>0$ with $(1-\e_i) \text{ sgn}(m_i)>0$
for all $i \in I_1$ such that $\a_i \bigl( \sum_{j \in I_1} (F^*)^r
\e_j^{r-1} m_j \o^\vee_j-w \cdot \mu \bigr)>0$ for all $i \in I_1$.
Now set $\l=\sum_{i \in I_1} \sum_{k=0}^{r-1} \e_i^k (F^*)^k m_i
\o^\vee_i$. Then $\a(\l)=\a(\mu)>0$ for all $\a>0$ with $w \a<0$.
Moreover, for $i \in I_1$ and $0 \le n<r$, \[\a_{\d^{-n} i} (F^*
\l-w \cdot \l)=\begin{cases} \a_{\d^{-n} i}(F^* \l-\l)=q^n m_i
\e_i^{n-1} (1-\e_i), & \text{ if } n \neq 0; \\ \a_i \bigl(
\sum_{j \in I_1} (F^*)^r \e_j^{r-1} m_j \o^\vee_j-w \cdot \mu
\bigr), & \text{ if } n=0. \end{cases}\]

Therefore, $\a_{\d^{-n} i} (F^*\l-w \cdot \l)>0$ for all $i \in I$ and $0 \le n<r$.

{\it Step 2}. Assume that $\Phi=\Phi_1 \sqcup \Phi_2$, where
$\Phi_k$ is generated by $I_k=I \cap \Phi_k$ and $\d(I_k)=I_k$ for
$k=1, 2$. Then $W=W_1 \times W_2$ and $\co=\co_1 \times \co_2$, where $\co_1$ (resp. $\co_2$) is a cuspidal
$\d$-conjugacy class of $W_1$ (resp. $W_2$). By assumption, there exists $w_k \in
(\co_k)_{\min}$ and $\mu_k=\sum_{i \in I_i} m_i \o^\vee_i$
for some $m_i \in \mathbb R$ such that $\a(\mu_k)>0$ for $\a>0$ with
$w_k \a<0$ and $\a_i (F^* \mu-w_k \cdot \mu)>0$ for all
$i \in I_k$. Now set $w=(w_1, w_2) \in \co_{\min}$ and
$\l=\mu_1+\mu_2$. Then $\a(\l)=\a(\mu_k)>0$ for $\a \in \Phi_k$ with
$\a>0$ and $w \a<0$. Also $\a_i \bigl(F^* \l-w \cdot \l)=\a_i
\bigl(F^* \mu_k-w_k \cdot \mu_k)>0$ for $i \in I_k$.

{\it Step 3}. Now we consider the general case. Here $\Phi=\sqcup \Phi^i$ and $\Phi^i=\sqcup \Phi^i_j$, where each $\Phi^i_j$ is irreducible and generated by $I^i_j=I \cap \Phi^i_j$ and for each $i$, $\d$ permutes $\{I^i_j\}$ cyclically. So we may apply Step 1 to
each $\Phi^i$ and then apply Step 2 to $\Phi$. It is easy to see that Lemma 1.11
holds in general.

\subsection*{1.13 Reduction to the condition $(J, w_1)$} It is easy to see that the map $w \mapsto w \i$ sends an $\d$-conjugacy class $\co$ to a $\d \i$-conjugacy class $\co^*$. If $\co$ is cuspidal, then so is $\co^*$. If $w \in \co_{\min}$, then $w \i \in \co^*_{\min}$. We will prove the following variant of Lemma 1.11.

(a). Let $\co$ be a cuspidal $\d \i$-conjugacy class. Then there exists $w \in \co_{\min}$ and $\mu \in X^\vee \otimes \mathbb R$ such that
$\a(\mu)>0$ for $\a>0$ with $w \i \a<0$ and $q \a_i(\mu)-(w \a_{\d \i(i)})(\mu)>0$ for all $i \in I$.

In fact, we will show that for most of the cases,

(b). there exists $w \in \co_{\min}$ and $\mu \in C$ such that $$q \a_i(\mu)-(w \a_{\d \i(i)})(\mu)>0, \qquad \text{ for all } i \in I.$$

The idea is as follows.

For $J \subset I$ and $x \in W^{\d \i(J)}$, set $$I(J, x, \d \i)=\max\{K \subset J; \Ad(x) \d \i(K)=K\}.$$ (In fact, if $K_1, K_2 \subset J$ with $\Ad(x) \d \i(K_i)=K_i$ for $i=1, 2$, then $\Ad(x) \d \i(K_1 \cup K_2)=K_1 \cup K_2$. Thus $\{K \subset J; \Ad(x) \d \i(K)=K\}$ contains a unique maximal element.)

By an observation in \cite[section 7]{He}, for each cuspidal
$\d \i$-conjugacy class $\co$ in a Weyl group of classical type, there
exists a maximal subset $J \subsetneqq I$, $w_1 \in W^{\d \i(J)}$ and a cuspidal $\Ad(w_1) \d \i$-conjugacy
class $\co'$ in $W_{I(J, w_1, \d \i)}$ such that for any $v \in
\co'_{\min}$, $v w_1 \in \co_{\min}$. We will see later in 1.15 that the observation is also valid for exceptional groups.

The following condition plays an essential role in our proof.

{\bf Condition $(J, w_1)$}: there exists $m_i \in \mathbb R_{>0}$ for $i \notin
I(J, w_1, \d \i)$ such that \[\tag{*} q \a_i (\sum_{j \notin I(J, w_1,
\d \i)} m_j \o^\vee_j)-(w_1 \a_{\d \i(i)}) (\sum_{j \notin I(J, w_1,
\d \i)} m_j \o^\vee_j)>0\] for $i \notin I(J, w_1, \d \i)$.

\begin{clm}
Keep notations as above. Suppose that the condition $(J, w_1)$ is true and that 1.13 (b) holds for $(\Phi_{I(J, w_1, \d \i)}, \Ad(w_1) \d \i, \co')$. Then 1.13(b) holds for $(\Phi_I, \d \i, \co)$.
\end{clm}

We simply write $K$ for $I(J, w_1, \d \i)$. By our assumption, there exists $v \in
\co'_{\min}$ and $m_i \in \mathbb R_{>0}$ for $i \in I$ such that for $i \in K$,
\begin{align*} & q \a_i(\sum_{j \in K} m_j \o^\vee_j)-(v \a_{w_1 \d \i (i)}) (\sum_{j \in K} m_j \o^\vee_j) \\ &=q \a_i(\sum_{j \in K} m_j \o^\vee_j)-(v w_1 \a_{\d \i (i)}) (\sum_{j \in K} m_j \o^\vee_j)>0.\end{align*}
and for $i \notin K$, \[q
\a_i (\sum_{j \notin K} m_j \o^\vee_j)-(w_1 \a_{\d \i(i)}) (\sum_{j
\notin K} m_j \o^\vee_j)>0.\]

Set $w=v w_1$ and $\l=\sum_{j \in K} m_j \o^\vee_j+m \sum_{j \notin K} m_j \o^\vee_j$ for $m \gg 0$. Notice that $w_1
\Phi_{\d \i(K)}=\Phi_K$. Thus for any $i \notin K$, $w_1
\a_{\d \i(i)} \notin \Phi_K$ and $w \a_{\d \i(i)}=v w_1 \a_{\d \i(i)}=w_1 \a_{\d \i(i)}+\sum_{j \in K}
a_j \a_j$ for $a_j \in \ZZ$ with $\sum_{j \in K} |a_j| \le n_0$.
Now for $i \in K$, \[q \a_i(\l)-(w \a_{\d \i(i)})(\l)=q \a_i(\sum_{j \in K} m_j \o^\vee_j)-(v w \a_{\d \i(i)}) (\sum_{j \in K} m_j
\o^\vee_j)>0\] and for $i \notin K$,
\begin{align*} & q \a_i (\l)-(w \a_{\d \i(i)})(\l) \\ &=m
\bigl(q \a_i (\sum_{j \notin K} m_j \o^\vee_j)-(w_1 \a_{\d \i(i)})
(\sum_{j \notin K} m_j \o^\vee_j) \bigr)-\sum_{j \in K} a_j m_j \\ & \ge m
\bigl(q \a_i (\sum_{j \notin K} m_j \o^\vee_j)-(w_1 \a_{\d \i(i)})
(\sum_{j \notin K} m_j \o^\vee_j) \bigr)-n_0 \max_{j \in K} m_j>0.\end{align*}

Therefore the 1.13(b) holds for $(\Phi_I, \d \i, \co)$.

\subsection*{1.14} We will show below that except the cases labelled with $\spadesuit$ (case 12 for type $E_8$, case 3 for type $F_4$, case 1 for type $G_2$, case 2 and case 4 for type $^2 F_4$), the condition $(J, w_1)$ are satisfied. Hence by induction on $|I|$, we can show that 1.13(b) holds for these cases. For the cases labelled with $\spadesuit$, we will prove 1.13(a) instead.

In \cite[Lemma 5.4 \& Lemma 5.7]{OR}, Orlik and Rapoport checked the condition $(J, w_1 \i)$ for type $A_n$ and $B_n$. In fact, in the case where $\d=id$, the condition $(J, w_1 \i)$ plays the same role in the proof of Lemma 1.11 as the condition $(J, w_1)$ does in the proof of 1.13(a). However, there are some big difference between the condition $(J, w_1)$ and the condition $(J, w_1 \i)$. In fact, there are more cases in the exceptional groups in which the condition $(J, w_1 \i)$ is not satisfied and for some of these cases, Lemma 1.11 is not easy to check directly. This is the reason why we prove 1.13(a) and the condition $(J, w_1)$ instead of Lemma 1.11 and the condition $(J, w_1 \i)$.

\subsection*{1.15} We use the same labelling of Dynkin diagram as in \cite{Bo}. We will use the same list of representatives of minimal length elements for all the cuspidal $\d \i$-conjugacy classes for the classical groups as in \cite[7.12-7.22]{He}. For the exceptional groups, we will also list a representative of minimal length elements for each cuspidal $\d \i$-conjugacy class. The representatives are presented as $v w_1$ for $w_1 \in W^{\d \i(J)}$ and $v$ is a minimal length element in the $\Ad(w_1) \d \i$-conjugacy class of $W_{I(J, w, \d \i)}$ that contains $v$. (These representatives are obtained by direct calculation based on the tables in \cite[Appendix B]{GP2} and \cite[section 6]{GKP}).

\ 

Set $$s_{[a, b]}=\begin{cases} s_a s_{a-1} \cdots s_b,
& \text{ if } a \ge b; \\ 0, & \text{ otherwise }. \end{cases}$$

\

{\bf Type $A_n$}

Set $J=I-\{1\}$. Here $w_1=s_{[n, 1]}$ and $I(J, w_1, \d \i)=\emptyset$. The inequality $(*)$ is just $q m_i-m_{i-1}>0$ for $i \neq 1$. So we may take $m_i=1$ for all $i$.

\

{\bf Type ${}^2 A_n$}

Set $J=I-\{n\}$. Here $w_1=s_{[n+1-a, 1]}$ for some $a \le \frac{n}{2}+1$ and $I(J,
w_1, \d \i)=\{a, a+1, \cdots, n-a\}$. The inequalities $(*)$ are
just $q m_i-m_{n+1-i}>0$ for $i<a-1$, $q m_{a-1}-m_{n+1-a}-m_{n+2-a}>0$ and $q m_i-m_{n-i}>0$ for $n-a<i<n$.
So we may take \[m_i=\begin{cases} 2, & \text{ if } i=a-1 \text{ or } n+1-a; \\ 1 & \text{ otherwise}. \end{cases}\]

\

{\bf Type $B_n$ and $C_n$}

Set $J=I-\{1\}$.

Case 1. $w_1=s_{[n-1, a]} \i s_{[n, 1]}$ for some $1 \le a<n$ and
$I(J, w_1, \d \i)=\{a+1, a+2, \cdots, n\}$. The inequalities $(*)$ are
just $q m_i-m_{i-1}>0$ for $1<i<a$ and $q m_a-m_{a-1}-m_a>0$. So we may take $m_i=1$ for $i<a$ and $m_a=2$.

Case 2. $w_1=s_{[n,1]}$ and $I(J, w_1, \d \i)=\emptyset$. The inequalities $(*)$ are just $q m_i-m_{i-1}>0$ for $1<i<n$ and $q m_n-m_n-\e m_{n-1}>0$, where \[\e=\begin{cases} 1, \quad & \text{ type } B_n;
\\ 2, \quad & \text{ type } C_n. \end{cases} \]

So we may take $m_n=3$ and $m_i=1$ for $i<n$.

\

{\bf Type $D_n$ and ${}^2 D_n$}

Set $J=I-\{\d(1)\}$.

Case 1. $w_1=s_{[n-2, a]} \i s_{[n, 1]}$ for some $a \le n-2$ and
$I(J, w_1, \d \i)=\{a+1, a+2, \cdots, n\}$. The inequalities $(*)$ are
just $q m_i-m_{i-1}>0$ for $1<i<a$ and $q m_a-m_{a-1}-m_a>0$. So we may take $m_a=1$ for $i<a$ and $m_a=2$.

Case 2. $w_1=s_{[n, 1]}$ and $I(J, w_1, \d \i)=\emptyset$. The
inequalities $(*)$ are just $q m_i-m_{i-1}>0$ for $1<i \le n-2$, $q m_{\d(n-1)}-m_{n-2}-m_n>0$ and $q
m_{\d(n)}-m_{n-2}-m_{n-1}>0$. So we may take $m_{n-1}=m_n=2$ and
$m_i=1$ for $i \le n-2$.

Case 3. $w_1=s_{[n-1, 1]}$ and $I(J, w_1, \d \i)=\emptyset$. The
inequalities $(*)$ are just $q m_i-m_{i-1}>0$ for $1<i \le n-2$, $q m_{\d(n-1)}-m_{n-2}>0$ and $q
m_{\d(n)}-m_{n-2}-m_{n-1}-m_n>0$. So we may take $m_{\d(n)}=3$ and $m_i=1$ for $i \neq \d(n)$.

\

{\bf Type ${}^3 D_4$}

Here $\d \i(s_1)=s_3$, $\d \i(s_3)=s_4$ and $\d \i(s_4)=s_1$. Set $J=I-\{4\}$.

Case 1. $w_1=s_2 s_1$ and $I(J, w_1, \d \i)=\emptyset$. The
inequalities $(*)$ are just $q m_1-m_2-m_3>0$, $q m_2-m_1>0$
and $q m_3-m_2-m_4>0$. So we may take $(m_1, m_2, m_3,
m_4)=(3,2,2,1)$.

Case 2. $w_1=s_{[3, 1]}$ and $I(J, w_1, \d \i)=\{1, 2\}$. The
inequality $(*)$ is just $q m_3-m_3-m_4>0$. So we may take $m_3=2$ and
$m_4=1$.

Case 3. $w_1=s_1 s_2 s_{[4, 1]}$ and $I(J, w_1, \d \i)=\{2, 3\}$. The
inequality $(*)$ is just $q m_1-m_4>0$. So we may take $m_1=m_4=1$.

\

{\bf Type $E_6$}

Set $J=I-\{6\}$.

Case 1. $w_1=s_{[6,1]} \i$ and $I(J, w_1, \d \i)=\emptyset$. The
inequalities $(*)$ are just $q m_1-m_3>0$, $q
m_2-m_1-m_3-m_4>0$, $q m_3-m_2-m_4>0$, $q m_4-m_5>0$ and $q
m_5-m_6>0$. So we may take $(m_1, m_2, m_3, m_4, m_5, m_6)=(2,4,3,1,
1, 1)$.

Case 2. $w_1=s_3 s_4 s_{[6,1]} \i$ and $I(J, w_1, \d \i)=\emptyset$..
The inequalities $(*)$ are just $q m_1-m_4>0$, $q
m_2-m_1>0$, $q m_3-m_2>0$, $q m_4-m_3-m_4-m_5>0$ and $q m_5-m_6>0$. So we
may take $(m_1, m_2, m_3, m_4, m_5, m_6)=(5,3,2,9,1,1)$.

Case 3. $w_1=s_2 s_4 s_5 s_3 s_4 s_{[6,1]} \i$, $I(J, w_1,
\d \i)=\{3,4\}$, $\Ad(w_1)$ is of order $2$ on $I(J, w_1, \d \i)$ and
$v=s_3$ or $s_3 s_4 s_3$. The inequalities $(*)$ are
just $q m_1-m_5>0$, $q m_2-m_1>0$ and $q m_5-m_1-m_5-m_6>0$. So we may
take $(m_1, m_2, m_5, m_6)=(3, 2, 5, 1)$.

Case 4. $w_1=w_0 w_0^J$, $I(J, w_1, \d \i)=\{2, 3, 4, 5\}$, $\Ad(w_1)$
is of order $3$ on $I(J, w_1, \d \i)$ and $\co'$ is cuspidal with $l(v)=8$. The inequality
$(*)$ is just $q m_1-m_6>0$. So we may take
$m_1=m_6=1$.

\

{\bf Type ${}^2 E_6$}

Set $J=I-\{1\}$.

Case 1. $w_1=s_2 s_{[6,4]} \i$ and $I(J, w_1, \d \i)=\emptyset$. The
inequalities $(*)$ are just $q m_2-m_4>0$, $q
m_3-m_6>0$, $q m_4-m_5>0$, $q m_5-m_2-m_3-m_4>0$ and $q m_6-m_1>0$.
So we may take $(m_1, m_2, m_3, m_4, m_5, m_6)=(1,2,1,3,5,1)$.

Case 2. $w_1=s_4 s_{[6,2]} \i$ and $I(J, w_1, \d \i)=\emptyset$. The
inequalities $(*)$ are just $q m_2-m_3>0$, $q
m_3-m_6>0$, $q m_4-m_4-m_5>0$, $q m_5-m_2>0$ and $q m_6-m_1-m_3-m_4>0$. So
we may take $(m_1, m_2, m_3, m_4, m_5, m_6)=(1,3,5,3,2,9)$.

Case 3. $w_1=s_5 s_4 s_{[6,2]} \i$, $I(J, w_1, \d \i)=\{4\}$ and $v=s_4$. The
inequalities $(*)$ are just $q m_2-m_3>0$, $q
m_3-m_5-m_6>0$, $q m_5-m_2>0$ and $q m_6-m_1-m_3-m_5>0$. So we may
take $(m_1, m_2, m_3, m_5, m_6)=(1,3,5,2,7)$.

Case 4. $w_1=s_{[6,4]} s_{[6,2]} \i$, $I(J, w_1, \d \i)=\{2, 3, 4,
5\}$, $\Ad(w_1) \d \i$ is of order $3$ on $I(J, w_1, \d \i)$, sending
$s_2$ to $s_3$, $s_3$ to $s_5$ and $s_5$ to $s_2$ and $\co'$ is cuspidal with $l(v)=4$ or
$6$. The inequality $(*)$ is just $q m_6-m_1-m_6>0$. So we
may take $m_1=1$ and $m_6=2$.

Case 5. $w_1=s_{[5,3]} s_{[6,4]} s_{[6,1]} \i$, $I(J, w_1, \d \i)=\{3,
4, 6\}$, $\Ad(w_1) \d \i$ is of order $2$ on $I(J, w_1, \d \i)$ and $v=s_3 s_4 s_3 s_6$. The inequalities $(*)$ is just $q m_2-m_1>0$. So we may take $(m_1, m_2, m_5)=(1,1,1)$.

Case 6. $w_1=s_3 s_1 w_0 w_0^J$, $I(J, w_1, \d \i)=\{5,6\}$, $\Ad(w_1) \d \i$ acts trivially on $I(J, w_1, \d \i)$ and $v=s_5 s_6$. The
inequalities $(*)$ are just $q m_2-m_1>0$, $q
m_3-m_2>0$ and $q m_4-m_3-m_4>0$. So we may take $(m_1, m_2, m_3,
m_4)=(1,1,1,2)$.

Case 7. $w_1=s_1 w_0 w_0^J$, $I(J, w_1, \d \i)=\{4,5,6\}$, $\Ad(w_1) \d \i$ acts trivially on $I(J, w_1, \d \i)$ and $v=s_{[6,4]}$. The
inequalities $(*)$ are just $q m_2-m_1-m_3>0$ and $q
m_3-m_2>0$. So we may take $m_1=1$ and $m_2=m_3=2$.

Case 8. $w_1=w_0 w_0^{\d \i(J)}$ and $v=w_0^{\d \i(J)}$. The inequalities $(*)$ are always satisfied.

\

{\bf Type $E_7$}

Set $J=I-\{7\}$.

Case 1. $w_1=s_{[7,1]} \i$ and $I(J, w_1, \d \i)=\emptyset$. The
inequalities $(*)$ are just $q m_1-m_3>0$, $q
m_2-m_1-m_3-m_4>0$, $q m_3-m_2-m_4>0$, $q m_4-m_5>0$, $q m_5-m_6>0$
and $q m_6-m_7>0$. So we may take $(m_1, m_2, m_3, m_4, m_5, m_6,
m_7)=(2,4,3,1,1,1,1)$.

Case 2. $w_1=s_3 s_4 s_{[7,1]} \i$ and $I(J, w_1, \d \i)=\emptyset$.
The inequalities $(*)$ are just $q m_1-m_4>0$, $q
m_2-m_1>0$, $q m_3-m_2>0$, $q m_4-m_3--m_4-m_5>0$, $q m_5-m_6>0$ and $q
m_6-m_7>0$. So we may take $(m_1, m_2, m_3, m_4, m_5, m_6,
m_7)=(5,3,2,9,1,1,1)$.

Case 3. $w_1=s_4 s_3 s_5 s_4 s_{[7,1]} \i$ and $I(J, w_1,
\d \i)=\emptyset$. The inequalities $(*)$ are just $q
m_1-m_5>0$, $q m_2-m_1>0$, $q m_3-m_2-m_4>0$, $2 m_4-m_3>0$,
$q m_5-m_4--m_5-m_6>0$ and $q m_6-m_7>0$. So we may take $(m_1, m_2, m_3,
m_4, m_5, m_6, m_7)=(5,3,3,2,4,1,1)$.

Case 4. $w_1=s_2 s_4 s_3 s_5 s_4 s_{[7,1]} \i$, $I(J, w_1, \d \i)=\{3,
4\}$, $\Ad(w_1)$ is of order $2$ on $I(J, w_1, \d \i)$ and $v=s_3$ or
$s_3 s_4 s_3$. The inequalities $(*)$ are just $q m_1-m_5>0$, $q m_2-m_1>0$, $q m_5-m_2-m_5-m_6>0$ and $q m_6-m_7>0$. So we
may take $(m_1, m_2, m_5, m_6, m_7)=(3,2,5,1,1)$.

Case 5. $w_1=s_3 s_4 s_2 s_{[5,3]} s_{[6,4]} s_{[7,1]} \i$, $I(J,
w_1, \d \i)=\{4\}$ and $v=s_4$. The inequalities $(*)$ are
just $q m_1-m_6>0$, $q m_2-m_1-m_3>0$, $q m_3-m_5>0$, $q m_5-m_2>0$
and $q m_6-m_3-m_5-m_6-m_7>0$. So we may take $(m_1, m_2, m_3, m_5, m_6,
m_7)=(7,5,2,3,7,1)$.

Case 6. $w_1=s_1 s_3 s_4 s_2 s_{[5,3]} s_{[6,4]} s_{[7,1]} \i$,
$I(J, w_1, \d \i)=\{2,3,4,5\}$, $\Ad(w_1)$ sends $\a_2$ to $\a_3$,
$\a_3$ to $\a_5$, $\a_4$ to $\a_4$, $\a_5$ to $\a_2$ and $\co'$ is cuspidal with $l(v)=4$,
$6$ or $8$. The inequalities $(*)$ are just $q
m_1-m_6>0$ and $q m_6-m_1-m_6-m_7>0$. So we may take $(m_1, m_6, m_7)=(3,5,1)$.

Case 7. $w_1=s_2 s_4 s_3 s_5 s_4 s_2 s_{[6,3]} s_{[7,4]} s_{[1,7]}
\i$, $I(J, w_1, \d \i)=\{3,4,5,6\}$, $\Ad(w_1)$ is of order $2$ on
$I(J, w_1, \d \i)$ and $v=w_0^{I(J, w_1, \d \i)}$. The inequalities $(*)$
are just $q m_1-m_7>0$ and $q m_2-m_1-m_2>0$. So we may take
$m_1=m_7=1$ and $m_2=2$.

Case 8. $w_1=s_{[6,4]} s_{[5,2]} \i s_1 s_3 s_4 s_2 s_{[6,3]}
s_{[7,4]} s_{[7,1]} \i$, $I(J, w_1, \d \i)=\{2,3,4,5\}$, $\Ad(w_1)$ is
of order $2$ on $I(J, w_1, \d \i)$ exchanging $\a_3$ and $\a_5$ and
$v=s_3 s_5 s_4 s_3 s_5 s_4 s_2$. The inequalities $(*)$
are just $q m_1-m_6-m_7>0$ and $q m_6-m_1>0$. So we may take $(m_1,
m_6, m_7)=(2,2,1)$.

Case 9. $w_1=w_0 w_0^J$ and $v=w_0^J$. The inequalities $(*)$ are
always satisfied.

\

{\bf Type $E_8$}

Set $J=I-\{8\}$.

Case 1. $w_1=s_{[8,1]} \i$ and $I(J, w_1, \d \i)=\emptyset$. The
inequalities $(*)$ are just $q m_1-m_3>0$, $q
m_2-m_1-m_3-m_4>0$, $q m_3-m_2-m_4>0$, $q m_4-m_5>0$, $q m_5-m_6>0$,
$q m_6-m_7>0$ and $q m_7-m_8>0$. So we may take $(m_1, m_2, m_3,
m_4, m_5, m_6, m_7, \break m_8)=(2,3,4,1,1,1,1,1)$.

Case 2. $w_1=s_3 s_4 s_{[8, 1]} \i$ and $I(J, w_1, \d \i)=\emptyset$.
The inequalities $(*)$ are just $q m_1-m_4>0$, $q m_2-m_1>0$, $q
m_3-m_2>0$, $q m_4-m_3-m_4-m_5>0$, $q m_5-m_6>0$, $q m_6-m_7>0$ and $q
m_7-m_8>0$. So we may take $(m_1, m_2, m_3, m_4, m_5, m_6, m_7,
m_8)=(5,3,2,9,1,1,1,1)$.

Case 3. $w_1=s_4 s_5 s_3 s_4 s_{[8, 1]} \i$ and $I(J, w_1,
\d \i)=\emptyset$. The inequalities $(*)$ are just $q
m_1-m_5>0$, $q m_2-m_1>0$, $q m_3-m_2-m_4>0$, $q m_4-m_3>0$,
$q m_5-m_4-m_5-m_6>0$, $q m_6-m_7>0$, $q m_7-m_8>0$. So we may take $(m_1,
m_2, m_3, m_4, m_5, m_6, m_7, m_8)=(5,3,3,2,4,1,1,1)$.

Case 4. $w_1=s_2 s_4 s_3 s_5 s_4 s_{[8,1]} \i$, $I(J, w_1, \d \i)=\{3,
4\}$, $\Ad(w_1)$ is of order $2$ on $I(J, w_1, \d \i)$ and $v=s_3$ or
$s_3 s_4 s_3$. The inequalities $(*)$ are just $q m_1-m_5>0$, $q m_2-m_1>0$, $q m_5-m_2-m_5-m_6>0$, $q m_6-m_7>0$ and $q m_7-m_8>0$. So we may take $(m_1, m_2, m_5, m_6, m_7,
m_8)=(3,2,5,1,1,1)$.

Case 5. $w_1=s_4 s_2 s_{[5,3]} s_{[6,4]} s_{[8, 1]} \i$ and $I(J,
w_1, \d \i)=\emptyset$. The inequalities $(*)$ are just $q
m_1-m_6>0$, $q m_2-m_1>0$, $q m_3-m_5>0$, $q m_4-m_3-m_4>0$, $q
m_5-m_2>0$, $q m_6-m_4-m_5-m_6-m_7>0$, $q m_7-m_8>0$. So we may take
$(m_1, m_2, m_3, m_4, m_5, m_6, m_7, m_8)=(9,5,2,3,3,17,1,1)$.

Case 6. $w_1=s_3 s_4 s_2 s_{[5,3]} s_{[6,4]} s_{[8, 1]} \i$, $I(J,
w_1, \d \i)=s_4$ and $v=s_4$. The inequalities $(*)$ are
just $q m_1-m_6>0$, $q m_2-m_1-m_3>0$, $q m_3-m_5>0$, $q m_5-m_2>0$,
$q m_6-m_3-m_5-m_6-m_7>0$ and $q m_7-m_8>0$. So we may take $(m_1, m_2, m_3,
m_5, m_6, m_7, m_8)=(7,5,2,3,13,1,1)$.

Case 7. $w_1=s_1 s_3 s_4 s_2 s_{[5,3]} s_{[6,4]} s_{[8, 1]} \i$,
$I(J, w_1, \d \i)=\{2,3,4,5\}$, $\Ad(w_1)$ sends $\a_2$ to $\a_3$,
$\a_3$ to $\a_5$, $\a_4$ to $\a_4$ and $\a_5$ to $\a_2$ and $\co'$ is cuspidal with
$l(v)=2$, $4$, $6$ or $8$. The inequalities $(*)$ are just $q m_1-m_6>0$,
$q m_6-m_1-m_6-m_7>0$ and $q m_7-m_8>0$. So we may take $(m_1, m_6, m_7,
m_8)=(3,5,1,1)$.

Case 8. $w_1=s_4 s_3 s_5 s_4 s_2 s_{[6,3]} s_{[7,4]} s_{[8, 1]} \i$,
$I(J, w_1, \d \i)=\{3, 6\}$, $\Ad(w_1)$ is of order $2$ on $I(J, w_1,
\d)$ and $v=s_3$. The inequalities $(*)$ are just $q
m_1-m_7>0$, $q m_2-m_1-m_4>0$, $q m_4-m_5>0$, $q m_5-m_2-m_4>0$ and
$q m_7-m_4-m_5-m_7-m_8>0$. So we may take $(m_1, m_2, m_4, m_5, m_7,
m_8)=(8,6,3,5,15,1)$.

Case 9. $w_1=s_2 s_4 s_3 s_5 s_4 s_2 s_{[6,3]} s_{[7,4]} s_{[8, 1]}
\i$, $I(J, w_1, \d \i)=\{3,4,5,6\}$, $\Ad(w_1)$ is of order $2$ on
$I(J, w_1, \d \i)$ and $v=s_3 s_4$ or $s_4 s_5 s_4 s_3$ or $w_0^{I(J,
w_1, \d)}$. The inequalities $(*)$ are just $q
m_1-m_7>0$, $q m_2-m_1-m_2>0$ and $q m_7-m_2-m_7-m_8>0$. So we may take $(m_1,
m_2, m_7, m_8)=(4,5,7,1)$.

Case 10. $w_1=s_5 s_4 s_{[7,2]} \i s_1 s_3 s_4 s_2 s_{[5,3]}
s_{[6,4]} s_{[8, 1]} \i$, $I(J, w_1, \d \i)=\{2, 4\}$, $\Ad(w_1)$ acts
trivially on $I(J, w_1, \d \i)$ and $v=s_2 s_4$. The inequalities $(*)$
are just $q m_1-m_7>0$, $q m_3-m_5-m_6>0$, $q
m_5-m_3>0$, $q m_6-m_1>0$, $q m_7-m_3-2 m_5-m_6-m_7-m_8>0$. So we may take
$(m_1, m_3, m_5, m_6, m_7, m_8)=(17,7,4,9,33,1)$.

Case 11. $w_1=s_{[6,1]} s_4 s_3 s_5 s_4 s_2 s_{[6,3]} s_{[7,4]}
s_{[8, 1]} \i$, $I(J, w_1, \d \i)=\{2,3,4,5\}$, $\Ad(w_1)$ is of order
$2$ on $I(J, w_1, \d \i)$ exchanging $\a_3$ and $\a_5$ and $v=s_2 s_4
s_5$ or $s_4 s_5 s_3 s_4 s_2 s_5 s_3$. The inequalities $(*)$ are just
$q m_1-m_6-m_7>0$, $q m_6-m_1>0$ and $q m_7-2 m_6-m_7-m_8>0$. So we may
take $(m_1, m_6, m_7, m_8)=(9,5,12,1)$.

Case 12  $\spadesuit$. $w_1=s_{[7,1]} s_4 s_3 s_5 s_4 s_2 s_{[6,3]} s_{[7,4]}
s_{[8, 1]} \i$, $I(J, w_1, \d \i)=I-\{7,8\}$, $\Ad(w_1)$ is of order
$2$ on $I(J, w_1, \d \i)$ and $\co'$ is cuspidal with $l(v)=12$, $14$, $16$, $18$ or $36$. In this case, the inequalities $(*)$ is never satisfied if $q=2$. However, notice that $w_1 \i v \i \a_8=w_1 \i \a_8>0$ for all $v \in W_{I(J, w_1, \d \i)}$. Thus if we choose $v$ to be a representative listed above in type ${}^2 E_6$ and $m_1, m_2, \cdots, m_6$ be the corresponding positive numbers there and take $m_7 \gg -m_8 \gg \max_{i=1, 2, \cdots, 6}\{m_i\}$, then one can see that 1.13(a) holds for $w=v w_1$ and $\mu=\sum_{i=1}^8 m_i \o^\vee_i$.

Case 13. $w_1=s_{[6,4]} s_{[6,2]} \i s_{[7,4]} s_{[6,2]} \i s_1 s_3
s_4 s_2 s_{[5,3]} s_{[8,4]} s_{[8,1]} \i$, $I(J, w_1,
\d \i)=\{2,3,4,5,7\}$, $\Ad(w_1)$ sending $\a_2$ to $\a_3$, $\a_3$ to
$\a_5$, $\a_4$ to $\a_4$, $\a_5$ to $\a_2$, $\a_7$ to $\a_7$ and $\co'$ is cuspidal with
$l(v)=9$. The inequalities $(*)$ are just $q
m_1-m_6-m_8>0$ and $q m_6-m_1-m_6>0$. So we may take $(m_1, m_6,
m_8)=(3,4,1)$.

Case 14. $w_1=s_3 s_4 s_2 s_{[7,5]} \i s_{[6,4]} \i s_{[5,3]} \i
s_{[3,1]} \i s_{[4,1]} s_{[5,3]} s_{[6,4]} s_2 s_{[7,3]} s_{[8,4]}
s_{[8,1]} \i$, $I(J, w_1, \d \i)=\{4,5,6,7\}$, $\Ad(w_1)$ acts
trivially on $I(J, w_1, \d \i)$ and $v=s_{[7,4]}$. The inequalities
$(*)$ are just $q m_1-m_3-m_8>0$, $q m_2-m_1-m_3>0$ and
$q m_3-m_2>0$. So we may take $(m_1, m_2, m_3, m_8)=(3,3,2,1)$.

Case 15. $w_1=s_1 s_3 s_4 s_2 s_{[7,5]} \i s_{[6,4]} \i s_3 s_4
s_{[4,1]} \i s_{[5,1]} s_4 s_3 s_{[6,4]} s_2 s_{[7,3]} s_{[8,4]}
s_{[8,1]} \i$, $I(J, w_1, \d \i)=\{2,3,4,5,6,7\}$, $\Ad(w_1)$ is of
order $2$ on $I(J, w_1, \d \i)$ and $v=s_3 s_4 s_{[5,2]} \i s_{[6,4]}
s_{[7,2]} \i$ or $s_{[4,2]} \i s_{[5,2]} \i s_4 s_{[5,2]} \i
s_{[6,4]} s_{[7,2]} \i$. The inequality $(*)$ is just
$q m_1-m_1-m_8>0$. So we may take $m_1=2$ and $m_8=1$.

Case 16. $w_1=s_{[7,1]} s_4 s_3 s_5 s_4 s_2 s_{[6,3]} s_{[7,4]}
s_{[7,1]} \i s_{[8,1]} s_4 s_3 s_5 s_4 s_2 s_{[6,3]} s_{[7,4]}
s_{[8,1]} \i$, $I(J, w_1, \d \i)=\{1,2,3,4,5,6\}$, $\Ad(w_1)$ acts
trivially on $I(J, w_1, \d \i)$ and $\co'$ is cuspidal with $l(v)=24$. The inequality $(*)$ is just $q m_7-m_7-m_8>0$. So we may take $m_7=2$ and $m_8=1$.

Case 17. $w_1=w_0 w_0^J$ and $v=w_0^J$. The inequalities $(*)$ are
always satisfied.

\

{\bf Type $F_4$}

Set $J=I-\{4\}$.

Case 1. $w_1=s_{[4, 1]}$ and $I(J, w_1, \d \i)=\emptyset$. The
inequalities $(*)$ are just $q m_2-m_1>0$,
$q m_3-m_2-m_3-m_4>0$ and $q m_4-m_3>0$. So we may take $(m_1, m_2, m_3,
m_4)=(1,1,5,3)$.

Case 2. $w_1=s_3 s_2 s_{[4, 1]}$ and $I(J, w_1, \d \i)=\emptyset$. The
inequalities $(*)$ are just $q m_2-m_1-m_2-2 m_3>0$, $q
m_3-m_4>0$ and $q m_4-m_2-m_3>0$. So we may take $(m_1, m_2, m_3,
m_4)=(1,12,5,9)$.

Case 3 $\spadesuit$. $w_1=s_2 s_3 s_2 s_{[4, 1]}$, $I(J, w_1, \d \i)=\{3, 4\}$, $\Ad(w_1)$ is of order $2$ on $I(J, w_1, \d \i)$ and $v=s_3$ or $s_3
s_4 s_3$. In this case, the inequalities $(*)$ is never satisfied if $q=2$. However, notice that $w_1 \i v \i \a_1=w_1 \i \a_1>0$ for all $v \in W_{I(J, w_1, \d \i)}$. Thus if we choose $v$ to be a representative listed above in type ${}^2 A_n$ and $m_3, m_4$ be the corresponding positive numbers there and take $m_2 \gg -m_1 \gg \max\{m_3, m_4\}$, then one can see that 1.13(a) holds for $w=v w_1$ and $\mu=\sum_{i=1}^4 m_i \o^\vee_i$.

Case 4. $w_1=s_{[3, 1]} s_3 s_2 s_{[4, 1]}$, $I(J, w_1, \d \i)=\{2\}$
and $v=s_2$. The inequalities $(*)$ are just $q m_3-m_3-m_4>0$ and $q m_4-m_1-m_3>0$. So we may take $(m_1, m_3, m_4)=(1,4,3)$.

Case 5. $w_1=s_{[4, 1]} s_3 s_2 s_{[4, 1]}$, $I(J, w_1, \d \i)=\{2,
3\}$, $\Ad(w_1)$ acts trivially on $I(J, w_1, \d \i)$ and $v=s_2
s_3$ or $s_2 s_3 s_2 s_3$. The inequality $(*)$ is just $q m_4-m_1-m_4>0$. So we may take $m_1=1$ and $m_4=2$.

Case 6. $w_1=s_1 w_0 w_0^{\d(J)}$, $I(J, w_1, \d \i)=\{3, 4\}$, $\Ad(w_1)$ acts trivially on $I(J, w_1, \d \i)$ and $v=s_3 s_4$. The
inequality $(*)$ is just $q m_2-m_1-m_2>0$. So we may take $m_1=1$ and
$m_2=2$.

Case 7. $w_1=w_0 w_0^J$ and $v=w_0^J$. The inequalities $(*)$ are always satisfied.

\

{\bf Type $G_2$}

Set $J=\{1\}$.

Case 1 $\spadesuit$. $w_1=s_1 s_2$ and $I(J, w_1, \d \i)=\emptyset$. In this case, the inequalities $(*)$ is never satisfied if $q=2$. However, take $m_1 \gg -m_2>0$, then 1.13(a) holds for $w_1$ and $\mu=m_1 \o^\vee_1+m_2 \o^\vee_2$.

Case 2. $w_1=s_1 s_2 s_1 s_2$ and $I(J, w_1, \d \i)=\emptyset$. The
inequality $(*)$ is just $q m_1-m_1-m_2>0$. So we may take $m_1=2$ and
$m_2=1$.

Case 3. $w_1=w_0$ and $I(J, w_1, \d \i)=\emptyset$. The inequalities
$(*)$ are always satisfied.

\

{\bf Type ${}^2 B_2$}

Set $J=\{1\}$.

Case 1. $w_1=s_1$ and $I(J, w_1, \d \i)=\emptyset$. The inequality
$(*)$ is just $q m_1-m_1-m_2>0$. So we may take $m_1=3$ and
$m_2=1$.

Case 2. $w_1=s_1 s_2 s_1$ and $I(J, w_1, \d \i)=\emptyset$. The
inequality $(*)$ is just $q m_1-m_2>0$. So we may take
$m_1=m_2=1$.

\

{\bf Type ${}^2 G_2$}

Set $J=\{2\}$.

Case 1. $w_1=s_2$ and $I(J, w_1, \d \i)=\emptyset$. The inequality
$(*)$ is just $q m_2-m_1-m_2>0$. So we may take $m_1=1$ and
$m_2=2$.

Case 2. $w_1=s_2 s_1 s_2$ and $I(J, w_1, \d \i)=\emptyset$. The inequality
$(*)$ is just $q m_2-2 m_1-m_2>0$. So we may take $m_1=1$ and
$m_2=3$.

Case 3. $w_1=s_2 s_1 s_2 s_1 s_2$ and $I(J, w_1, \d \i)=\emptyset$. The inequality
$(*)$ is just $q m_2-m_1>0$. So we may take $m_1=m_2=1$.

\

{\bf Type ${}^2 F_4$}

Set $J=I-\{4\}$.

Case 1. $w_1=s_2 s_1$ and $I(J, w_1, \d \i)=\emptyset$. The inequality
$(*)$ are just $q m_1-m_4>0$, $q m_2-m_2-m_3>0$ and $q m_3-m_1>0$. So we may take $(m_1, m_2, m_3, m_4)=(1,3,1,1)$.

Case 2  $\spadesuit$. $w_1=s_2 s_{[3,1]}$ and $I(J, w_1, \d \i)=\emptyset$. In this case, the inequalities $(*)$ is never satisfied if $q=\sqrt{2}$. However, take $-m_4 \gg m_2=m_3 \gg m_1>0$, then 1.13(a) holds for $w_1$ and $\mu=\sum_{i=1}^4 m_i \o^\vee_i$.

Case 3. $w_1=s_1 s_2 s_{[3,1]}$, $I(J, w_1, \d \i)=\{2, 3\}$, $\Ad(w_1) \d \i$ is of order $2$ on $I(J, w_1, \d \i)$ and $v=s_2$ or $s_2
s_3 s_2$. The inequality $(*)$ is just $q m_1-m_1-m_4>0$. So
we may take $m_1=3$ and $m_4=1$.

Case 4 $\spadesuit$. $w_1=s_{[3,1]} s_2 s_3 s_2 s_{[4,1]}$ and $I(J, w_1,
\d \i)=\emptyset$. In this case, the inequalities $(*)$ is never satisfied if $q=\sqrt{2}$. However, 1.13(a) holds for $w=w_1$ and $\mu=3 \o^\vee_1+\o^\vee_2+3 \o^\vee_3-3 \o^\vee_4$.

Case 5. $w_1=s_2 s_{[3,1]} s_2 s_3 s_2 s_{[4,1]}$, $I(J, w_1,
\d)=\{1,3\}$, $\Ad(w_1) \d \i$ is of order $2$ on $I(J, w_1, \d \i)$
and $v=s_2$. The inequality $(*)$ is just $q m_2-m_2-m_4>0$. So we may take $m_2=3$ and $m_4=1$.

Case 6. $w_1=w_0 w_0^{\d \i(J)}$, $I(J, w_1, \d \i)=\{2, 3\}$, $\Ad(w_1) \d \i$ is of order $2$ on $W_{I(J, w_1, \d \i)}$ and $v=s_2 s_3
s_2$. The inequality $(*)$ is just $q m_1-m_4>0$. So we may
take $m_1=m_4=1$.

\section*{Acknowledgements}
We thank G. Lusztig for helpful discussions on Deligne-Lusztig varieties. We thank S. Orlik and M. Rapoport for many useful comments and suggestions for improvement in an earlier version of this paper. We also thank J. Starr for helpful discussions on affine schemes. The computations for exceptional groups were done with the aid of CHEVIE \cite{CH} package of GAP \cite{GA}.

\bibliographystyle{amsalpha}

\end{document}